\documentclass[a4paper,reqno]{amsbook}%
\usepackage{amsfonts}
\usepackage{amsmath}
\usepackage{amssymb}
\usepackage[a4paper,asymmetric,includehead,nofoot,verbose]{geometry}%
\usepackage{graphicx}
\theoremstyle{plain}

\numberwithin{equation}{section}
\begin{document}
\frontmatter
\title[Overview]{Overview document for:\\A weight function theory of basis function interpolants and smoothers.}
\author{Phillip Y. Williams}
\address{5/31 Moulden Court, Belconnen, ACT 2617.}
\email{philwill@pcug.org.au, PhilWilliams@uc.edu.au}
\thanks{I would like to thank the University of Canberra for their generous assistance.}
\thanks{.}
\thanks{Also, thanks to my Masters degree supervisors Dr Markus Hegland and Dr Steve
Roberts of the CMA at the Australian National University.}
\date{October 2007}
\tableofcontents

\begin{abstract}
This document is a brief overview of two documents which continue to develop
the weight function theory of basis function smoothers and interpolants. One
document considers the zero order theory and one considers the positive order theory.

\end{abstract}
\maketitle

\subsection{Change register}

\ \medskip

\textbf{05/Jul/2007}\quad Created this document using the abstracts from the
version 1 documents.

\textbf{22/Oct/2007}\quad Altered this document using the abstracts from the
consolidated version 2 documents i.e. version 2 of arXiv:0708.0780 and version
2 of arXiv:0708.0795.

\mainmatter

\section{Overview}

This work currently consists of \textbf{two documents} which continue to
develop the \textbf{Light weight function theory of basis function smoothers
and interpolants}. One document considers the \textbf{zero order theory} and
one considers the \textbf{positive order theory}.

In brief, some important general features:\smallskip

\begin{enumerate}
\item Extends the positive order work of Light and Wayne
\cite{LightWayne98PowFunc}, \cite{LightWayneX98Weight} and
\cite{LightWayne95ErrEst} to the zero order case and extends the positive
order case to tensor product weight functions.\medskip

For both the positive and zero order cases:\smallskip

\item A weight function is first defined and then used to define a continuous
basis function and a data function Hilbert space are defined using the Fourier
transform. This technique is illustrated by several examples.

\item The standard minimal norm and seminorm interpolants are defined and
pointwise orders of convergence are derived on a bounded set.

\item We define the well known variational non-parametric smoother which
stabilizes the interpolant using a smoothing parameter - I call this the Exact
smoother. Orders of uniform pointwise convergence are derived on a bounded
open set.

\item A scalable smoother is derived which I call the Approximate smoother.
Orders of uniform pointwise convergence are derived on a bounded open
set.\medskip

\item For the \textbf{zero order} case numeric examples are given which
compare the theoretical and actual errors w.r.t. the data function.
\end{enumerate}

\section{Zero order document (arXiv:0708.0780)}

Here is a short description of the document (the abstract).\medskip

In this document I develop a weight function theory of zero order basis
function interpolants and smoothers.\smallskip

In \textbf{Chapter 1} the basis functions and data spaces are defined directly
using weight functions. The data spaces are used to formulate the variational
problems which define the interpolants and smoothers discussed in later
chapters. The theory is illustrated using some standard examples of radial
basis functions and a class of weight functions I will call the tensor product
extended B-splines.

In \textbf{Chapter 2} the theory of Chapter 1 is used to prove the pointwise
convergence of the minimal norm basis function interpolant to its data
function and to obtain orders of convergence. The data functions are
characterized locally as Sobolev-like spaces and the results of several
numerical experiments using the extended B-splines are presented.

In \textbf{Chapter 3} a large class of tensor product weight functions will be
introduced which I call the central difference weight functions. These weight
functions are closely related to the extended B-splines and have similar
properties. The theory of this document is then applied to these weight
functions to obtain interpolation convergence results. To understand the
theory of interpolation and smoothing it is not necessary to read this chapter.

In \textbf{Chapter 4} a non-parametric variational smoothing problem will be
studied using the theory of this document with special interest in its order
of pointwise convergence of the smoother to its data function. This smoothing
problem is the minimal norm interpolation problem stabilized by a smoothing coefficient.

In \textbf{Chapter 5} a non-parametric, scalable, variational smoothing
problem will be studied, again with special interest in its order of pointwise
convergence to its data function. We discuss the \textit{SmoothOperator}
software (freeware) package which implements the Approximate smoother
algorithm. It has a full user manual which describe several tutorials and data experiments.

\section{Positive order document (arXiv:0708.0795)}

Here is a short description of the document (the abstract).\medskip

In this document I develop a weight function theory of positive order basis
function interpolants and smoothers.

In \textbf{Chapter 1} the basis functions and data spaces are defined directly
using weight functions. The data spaces are used to formulate the variational
problems which define the interpolants and smoothers discussed in later
chapters. The theory is illustrated using some standard examples of radial
basis functions and a class of weight functions I will call the tensor product
extended B-splines.

\textbf{Chapter 2} shows how to prove functions are basis functions without
using the awkward space of test functions So,n which are infinitely smooth
functions of rapid decrease with several zero-valued derivatives at the
origin. Worked examples include several classes of well-known radial basis functions.

The goal of \textbf{Chapter 3} is to derive `modified' inverse-Fourier
transform formulas for the basis functions and the data functions and to use
these formulas to obtain bounds for the rates of increase of these functions
and their derivatives near infinity.

In \textbf{Chapter 4} we prove the existence and uniqueness of a solution to
the minimal seminorm interpolation problem. We then derive orders for the
pointwise convergence of the interpolant to its data function as the density
of the data increases.

In \textbf{Chapter 5} a well-known non-parametric variational smoothing
problem will be studied with special interest in the order of pointwise
convergence of the smoother to its data function. This smoothing problem is
the minimal norm interpolation problem stabilized by a smoothing coefficient.

In \textbf{Chapter 6} a non-parametric, scalable, variational smoothing
problem will be studied, again with special interest in its order of pointwise
convergence to its data function.

\bibliographystyle{amsplain}
\bibliography{BasisFunc}

\providecommand{\bysame}{\leavevmode\hbox to3em{\hrulefill}\thinspace}
\providecommand{\MR}{\relax\ifhmode\unskip\space\fi MR }
\providecommand{\MRhref}[2]{%
  \href{http://www.ams.org/mathscinet-getitem?mr=#1}{#2}
}
\providecommand{\href}[2]{#2}
\begin{thebibliography}{1}

\bibitem{LightWayne95ErrEst}
W.~Light and H.~Wayne, \emph{Error estimates for approximation by radial basis
  functions.}, Approximation theory, wavelets and applications (Maratea, 1994),
  NATO Adv. Sci. Inst. Ser. C Math. Phys. Sci., vol. 454, Kluwer Acad. Publ.,
  Dordrecht, 1995, pp.~215--246.

\bibitem{LightWayne98PowFunc}
\bysame, \emph{On power functions and error estimates for radial basis function
  interpolation}, J. Approx. Th. \textbf{92} (1998), no.~2, 245--266.

\bibitem{LightWayneX98Weight}
\bysame, \emph{Spaces of distributions, interpolation by translates of a basis
  function and error estimates}, Numer. Math. \textbf{81} (1999), no.~3,
  415--450.

\end{thebibliography}

\end{document}